\documentstyle[11pt]{article}
  \newlength{\standardunitlength}
\newtheorem{cor}{Corollary} \newtheorem{lemma}{Lemma}
\newtheorem{theorem}{Theorem} 
\newenvironment{proof}{\noindent {\sc Proof:}}{$\Box$ \vspace{2 ex}}

\begin{document}

\begin{center} {\bf Finite Affine Groups: Cycle Indices,
Hall-Littlewood Polynomials, and Probabilistic Algorithms}
\end{center}

\begin{center}
By Jason Fulman
\end{center}

\begin{center}
Affiliation at time of writing: Stanford University
\end{center}

\begin{center}
Current affiliation: University of Pittsburgh
\end{center}

\begin{center}
Department of Mathematics
\end{center}

\begin{center}
301 Thackeray Hall
\end{center}

\begin{center}
Pittsburgh, PA, 15260
\end{center}

\begin{center}
fulman@math.pitt.edu
\end{center}

\begin{center}
Submitted: September 19, 2000
\end{center}

\begin{center}
Final revised version for J. Algebra: August 4, 2001
\end{center}

\newpage

\begin{center}
Proposed running head: Finite Affine Groups
\end{center}

\begin{center}
Please send proofs to:
\end{center}

\begin{center}
Jason Fulman
\end{center}

\begin{center}
University of Pittsburgh
\end{center}

\begin{center}
Department of Mathematics
\end{center}

\begin{center}
301 Thackeray Hall
\end{center}

\begin{center}
Pittsburgh, PA 15260
\end{center}

\begin{center}
fulman@math.pitt.edu
\end{center}

\newpage

\begin{abstract} The study of asymptotic properties of the conjugacy
class of a random element of the finite affine group leads one to
define a probability measure on the set of all partitions of all
positive integers. Four different probabilistic understandings of this
measure are given--three using symmetric function theory and one using
Markov chains. This leads to non-trivial enumerative results. Cycle
index generating functions are derived and are used to compute the
large dimension limiting probabilities that an element of the affine
group is separable, cyclic, or semisimple and to study the convergence
to these limits. The semisimple limit involves both Rogers-Ramanujan
identities. This yields the first examples of such computations for a
maximal parabolic subgroup of a finite classical group. \end{abstract}

\begin{center} Key words: Conjugacy class, classical group, affine
group, Hall-Littlewood polynomial, symmetric function. \end{center}

\section{Introduction} \label{Introduction}

	The conjugacy classes of the unitary groups with complex
entries are simply the set of its eigenvalues. Thus the enormous body
of recent work on eigenvalues of random matrices is to a large extent
a study of conjugacy classes. Given the power of random matrix models
(see for instance \cite{KS1} and \cite{KS2} which show the predictive
power of random matrices for the study of local properties of the
zeros of the Riemann zeta function), it is natural to study conjugacy
in random matrices over finite fields as well.

	The papers \cite{F1}, \cite{F2}, \cite{B}, \cite{F3} give a
purely probabilistic understanding of conjugacy classes in the finite
classical groups and in the group of upper triangular matrices over a
finite field. One outcome was a simple and motivated proof of the
Rogers-Ramanujan identities \cite{F4}.

	The current paper considers the affine group $A(n,q)$ (all
matrices in $GL(n+1,q)$ with $x_{11}=1$ and $x_{j1}=0$ for $j \geq 2$)
and the maximal parabolic subgroup $P(n,q)$ (all matrices in
$GL(n+1,q)$ with $x_{11} \neq 0$ and $x_{j1}=0$ for $j \geq 2$). 

	Section \ref{S2} derives an expression for the chance that an
element of $A(n,q)$ or $P(n,q)$ has a given rational canonical form in
$GL(n+1,q)$, giving rise to cycle index generating functions for these
groups. Section \ref{S2} also makes useful contact with the
Hall-Littlewood symmetric functions.

	Section \ref{S3} opens by defining a new and natural
probability measure $N_{u,q}$ on the set of all partitions of all
positive integers and by relating it to an analogous measure for the
finite general linear groups. Then it gives four purely probabilistic
algorithms for ``growing'' random partitions according to $N_{u,q}$
and also an algorithm which arises when considering conjugacy classes
only of unipotent elements in $A(n,q)$ or $P(n,q)$.

	Section \ref{S4} gives non-trivial enumerative applications of
the algorithms of Section \ref{S3}. For instance the $n \rightarrow
\infty$ chance that an element of $A(n,q)$ or $P(n,q)$ has a fixed
space of a given dimension is shown to have a simple product formula,
where each factor in the product has a clear probabilistic
meaning. The number of unipotent elements in $A(n,q)$ or $P(n,q)$ of a
given rank is computed.

	Section \ref{S5} of this paper uses the generating functions
of Section \ref{S2} to calculate the probabilities that an element of
$A(n,q)$ or $P(n,q)$ is separable, cyclic, or semisimple. Analogous
results are known for the finite classical groups
\cite{F5},\cite{W},\cite{FNP} and are important for computational
group theory \cite{NP},\cite{NP4}. The motivation for Section \ref{S5}
is the fact that probabilistic estimates in maximal subgroups of the
finite classical groups are also useful for group theory. To state a
typical result, recall that a matrix is called cyclic if its minimal
polynomial is equal to its characteristic polynomial. Whereas the $n
\rightarrow \infty$ limiting probability that an element of $GL(n,q)$
is cyclic is equal to $\frac{1-1/q^5}{1+1/q^3}$, Section 5 shows that
the corresponding probability for $A(n,q)$ or $P(n,q)$ is
$\frac{1-1/q}{1-1/q+1/q^2} \frac{1-1/q^5}{1+1/q^3}$, which is
asymptotically $1-1/q^2+O(1/q^3)$ rather than $1-1/q^3+O(1/q^4)$. It
is worth emphasizing that at present even for $GL(n,q)$ the only
method for performing such exact calculations is through the use of
generating functions.

\section{Cycle Indices, Partitions and Hall-Littlewood Polynomials}
\label{S2}

	To begin we recall some standard notation about partitions
which will be used throughout the paper. Let $\lambda$ be a partition
of some non-negative integer $|\lambda|$ into parts $\lambda_1 \geq
\lambda_2 \geq \cdots$. Let $m_i(\lambda)$ be the number of parts of
$\lambda$ of size $i$, and let $\lambda'$ be the partition dual to
$\lambda$ in the sense that $\lambda_i' = m_i(\lambda) +
m_{i+1}(\lambda) + \cdots$. Let $n(\lambda)$ be the quantity $\sum_{i
\geq 1} (i-1) \lambda_i$. It is also useful to define the diagram
associated to $\lambda$ as having the $j$th row consist of $\lambda_j$
boxes. For example the diagram of the partition $(5441)$ is:

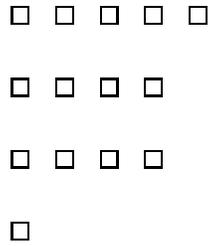
\[ \begin{array}{c c c c c}
		\framebox{} & \framebox{} & \framebox{} & \framebox{} & \framebox{}  \\
		\framebox{} & \framebox{} & \framebox{} & \framebox{} &    \\
		\framebox{} & \framebox{} & \framebox{} & \framebox{} &    \\
		\framebox{} & & & &  
	  \end{array} \] The notation $P_{\lambda}(x_1,x_2,\cdots;t)$
denotes the Hall-Littlewood symmetric function. The reader is referred
to Chapter 3 of \cite{Mac} for a thorough discussion of its
properties. The symbol $(\frac{1}{q})_i$ will denote $(1-1/q) \cdots
(1-1/q^i)$.

	Next recall (e.g. Chapter 6 of \cite{H}) that the conjugacy
classes of $GL(n,q)$ are parameterized by rational canonical form. This
form corresponds to the following combinatorial data. To each monic
non-constant irreducible polynomial $\phi$ over the field of $q$
elements $F_q$, associate a partition (perhaps the trivial partition)
$\lambda_{\phi}$ of some non-negative integer $|\lambda_{\phi}|$. Let
$deg(\phi)$ denote the degree of $\phi$. The only restrictions
necessary for this data to represent a conjugacy class are that the
partition corresponding to the polynomial $\phi(z)=z$ is empty and that
$\sum_{\phi} |\lambda_{\phi}| deg(\phi) = n.$

	This conjugacy data arises explicitly as follows. Given an element $\alpha \in GL(n,q)$ operating on the vector space $V$, there is a unique direct sum decomposition $V= \bigoplus V_{\phi}$ where the characteristic polynomial of $\alpha$ on $V_{\phi}$ is a power of $\phi$ and the characteristic polynomials on any two summands are relatively prime. Furthermore $\alpha$ decomposes each $V_{\phi}$ into a sum of cyclic subspaces. This decomposition need not be unique but the dimensions of the cyclic subspaces in the decomposition of $V_{\phi}$ are unique and are the row lengths of the partition $\lambda_{\phi}$. For example, the identity matrix has $\lambda_{z-1}$ equal to $(1^n)$ and an elementary reflection with $a
\neq 0$ in the $(1,2)$ position, ones on the diagonal and zeros
elsewhere has $\lambda_{z-1}$ equal to $(2,1^{n-2})$. As another
example, the characteristic polynomial of an element with conjugacy
class data $\{\lambda_{\phi}\}$ is $\prod_{\phi \neq z}
\phi^{|\lambda_{\phi}|}$.

	The first aim of this section is to find expressions for the
chance that an element of $A(n,q)$ or $P(n,q)$ has given rational
canonical form data. This is a cruder invariant than conjugacy in
$A(n,q)$ or $P(n,q)$. However most conjugacy class functions on
$A(n,q)$ or $P(n,q)$ of interest depend only on the Jordan form in
$GL(n+1,q)$. (A similar reduction proved useful \cite{B},\cite{F3} for
upper triangular matrices over a finite field, for which the theory of
wild quivers implies that there is no finite parameterization of
conjugacy classes).

	To find the chance that an element of $A(n,q)$ or $P(n,q)$ has
a given rational canonical form we use a result of Nakada and Shinoda
\cite{NaS} on the parameterization and sizes of conjugacy classes in
$A(n,q)$ or $P(n,q)$. Then a combinatorial reformulation followed
by a summation gives the result we seek.

	For the statement of Theorem \ref{classes}, we use the
notation that $|Z_G(c)|$ is the size of the centralizer in the group
$G$ of any element in the conjugacy class $c$. The symbol $F_q^*$
denotes the non-zero elements in $F_q$. We use the convention that $GL(0,q)$ consists of one conjugacy class of size one.

\begin{theorem} \label{classes} (\cite{NaS}) 

\begin{enumerate}

\item The conjugacy classes of $A(n,q)$ are parameterized by pairs
$(c_{n+1-k},k)$ where $0< k \leq n+1$ and $c_{n+1-k}$ is a conjugacy
class of $GL(n+1-k,q)$. Letting $\{\overline{\lambda_{\phi}}\}$ be the
conjugacy class data corresponding to $c_{n+1-k}$, the $GL(n+1,q)$
rational canonical form of the corresponding class in $A(n,q)$ is
given by $\lambda_{\phi}=\overline{\lambda_{\phi}}$ for $\phi \neq z-1$ and
by letting $\lambda_{z-1}=(k) \cup \overline{\lambda_{z-1}}$ be the
partition formed by adding a row of length $k$ to
$\overline{\lambda_{z-1}}$. The centralizer size of the corresponding
conjugacy class is

\[ |Z_{GL(n+1-k,q)}(c_{n+1-k})| q^{k-1+2 \sum_{i=1}^{k-1} i
m_i(\overline{\lambda_{z-1}}) + (2k-1) \sum_{i=k}^{n-k}
m_i(\overline{\lambda_{z-1}})} .\]

\item The conjugacy classes of $P(n,q)$ are parameterized by triples
$(c_{n+1-k},k,a)$ where $0< k \leq n+1$, $a \in F_q^*$ and $c_{n+1-k}$ is
a conjugacy class of $GL(n+1-k,q)$. Letting $\{\overline{\lambda_{\phi}}\}$
be the conjugacy class data corresponding to $c_{n+1-k}$, the
$GL(n+1,q)$ rational canonical form of the corresponding class in
$A(n,q)$ is given by $\lambda_{\phi}=\overline{\lambda_{\phi}}$ for $\phi
\neq z-a$ and by letting $\lambda_{z-a}=(k) \cup \overline{\lambda_{z-a}}$
be the partition formed by adding a row of length $k$ to
$\overline{\lambda_{z-a}}$. The centralizer size of the corresponding
conjugacy class is

\[ |Z_{GL(n+1-k,q)}(c_{n+1-k})| (q-1) q^{k-1+2 \sum_{i=1}^{k-1} i
m_i(\overline{\lambda_{z-1}}) + (2k-1) \sum_{i=k}^{n-k}
m_i(\overline{\lambda_{z-1}})} .\]
\end{enumerate}

\end{theorem}

	Lemma \ref{useful} below was Theorem 9 in the thesis \cite{F0}. Of
course the formula for the conjugacy class sizes in $GL(n,q)$ is not
due to the author (see for instance page 219 of \cite{Mac} for the
third expression), but the combinatorial rewritings of it in Lemma
\ref{useful} are very useful.

\begin{lemma} \label{useful} (\cite{F0}) The conjugacy class of
$GL(n,q)$ corresponding to the data $\{\lambda_{\phi}\}$ has size

\[ \frac{|GL(n,q)|}{\prod_{\phi \neq z} c_{GL,\phi,q}(\lambda_{\phi})} \] where
\begin{eqnarray*}
c_{GL,\phi,q}(\lambda) & = & q^{2deg(\phi) [\sum_{h<i} h
m_h(\lambda) m_i(\lambda) + \frac{1}{2} \sum_i (i-1) m_i(\lambda)^2]}
\prod_i |GL(m_i(\lambda),q^{deg(\phi)})|\\
& = & q^{deg(\phi) [\sum_i
(\lambda'_i)^2]} \prod_{i} (\frac{1}{q^{deg(\phi)}})_{m_i(\lambda)}\\
& = &
\frac{q^{deg(\phi)
n(\lambda)}}{P_{\lambda}(\frac{1}{q^{deg(\phi)}},\frac{1}{q^{2deg(\phi)}},\cdots;\frac{1}{q^{deg(\phi)}})}
\end{eqnarray*}

\end{lemma}

	Theorem \ref{probabform} is our first main result. The
notation used in the proof is the same as in the statement of Theorem
\ref{classes}. Recall that $\lambda_{\phi,1}'$ is the size of the first column of $\lambda_{\phi}$.

\begin{theorem} \label{probabform}

\begin{enumerate}

\item The number of elements of $A(n,q)$ with $GL(n+1,q)$ rational
canonical form data $\{\lambda_{\phi}\}$ is equal to

\[ \frac{|A(n,q)| (q^{\lambda_{z-1,1}'}-1)} {\prod_{\phi \neq z}
q^{deg(\phi) \cdot \sum_i (\lambda_{\phi,i}')^2}
(\frac{1}{q^{deg(\phi)}})_{m_i(\lambda_{\phi})}}.\] 

\item The number of elements of $P(n,q)$ with $GL(n+1,q)$ rational
canonical form data $\{\lambda_{\phi}\}$ is equal to

\[ \sum_{a \in F_q^*} \frac{|P(n,q)| (q^{\lambda_{z-a,1}'}-1)}{q-1}
\frac{1}{\prod_{\phi \neq z} q^{deg(\phi) \cdot \sum_i (\lambda_{\phi,i}')^2}
(\frac{1}{q^{deg(\phi)}})_{m_i(\lambda_{\phi})}}.\]

\end{enumerate}
\end{theorem}

\begin{proof} First we consider $A(n,q)$. From part 1 of Theorem
\ref{classes} the number of elements of $A(n,q)$ with $GL(n+1,q)$
rational canonical form data $\{\lambda_{\phi}\}$ is 

\[ \sum_{(c_{n+1-k},k) \atop k \geq 1, (k) \cup \overline{\lambda_{z-1}} =
\lambda_{z-1}} \frac{|A(n,q)|}{|Z_{GL(n+1-k,q)}(c_{n+1-k})| q^{k-1+2
\sum_{i=1}^{k-1} i m_i(\overline{\lambda_{z-1}}) + (2k-1) \sum_{i=k}^{n-k}
m_i(\overline{\lambda_{z-1}})}}.\] From Lemma \ref{useful},
$|Z_{GL(n+1-k,q)}(c_{n+1-k})|$ can be written as \[\prod_{\phi \neq z}
\prod_i q^{deg(\phi) \cdot (\overline{\lambda_{\phi,i}'})^2}
(\frac{1}{q^{deg(\phi)}})_{m_i(\overline{\lambda_{\phi}})}.\] Thus the
number of elements of $A(n,q)$ with $GL(n+1,q)$ rational canonical
form data $\{\lambda_{\phi}\}$ is

\begin{eqnarray*}
|A(n,q)| \sum_{(k,\overline{\lambda_{z-1}}) \atop k \geq 1,(k) \cup
\overline{\lambda_{z-1}}=\lambda_{z-1}} \frac{1}{q^{k-1+2 \sum_{i=1}^{k-1}
i m_i(\overline{\lambda_{z-1}}) + (2k-1) \sum_{i=k}^{n-k}
m_i(\overline{\lambda_{z-1}})}} \frac{1}{\prod_i q^{ (\overline{\lambda_{z-1,i}'})^2}
(\frac{1}{q})_{m_i(\overline{\lambda_{z-1}})} }\\
\cdot \prod_{\phi \neq z,z-1} \prod_i \frac{1}{q^{deg(\phi) \cdot
(\overline{\lambda_{\phi,i}'})^2}
(\frac{1}{q^{deg(\phi)}})_{m_i(\overline{\lambda_{\phi}})}}
\end{eqnarray*}

	Consequently it is sufficient to prove that

\begin{eqnarray*}
& & \sum_{(k,\overline{\lambda_{z-1}}) \atop k \geq 1,(k) \cup
\overline{\lambda_{z-1}}=\lambda_{z-1}} \frac{1}{q^{k-1+2 \sum_{i=1}^{k-1}
i m_i(\overline{\lambda_{z-1}}) + (2k-1) \sum_{i=k}^{n-k}
m_i(\overline{\lambda_{z-1}})} \prod_i q^{ (\overline{\lambda_{z-1,i}'})^2}
(\frac{1}{q})_{m_i(\overline{\lambda_{z-1}})}}\\
& = & \frac{(q^{\lambda_{z-1,1}'}-1)} {\prod_{i} q^{
(\lambda_{z-1,i}')^2} (\frac{1}{q})_{m_i(\lambda_{z-1})}}.
\end{eqnarray*} Since
only partitions corresponding to the polynomial $z-1$ are involved in
this last equation, we simplify notation by suppressing the
dependence on $z-1$ and by further replacing $m_i(\lambda)$ and
$m_i(\overline{\lambda})$ by $m_i$ and $\overline{m_i}$.

	From Lemma \ref{useful} one sees that

\begin{eqnarray*}
& & \sum_{(k,\overline{\lambda}) \atop k \geq 1,(k) \cup \overline{\lambda}=\lambda} \frac{1}{q^{k-1+2
\sum_{i=1}^{k-1} i \overline{m_i} + (2k-1) \sum_{i=k}^{n-k} \overline{m_i}} \prod_i
q^{ (\overline{\lambda}_i')^2} (\frac{1}{q})_{\overline{m_i}}}\\
& = & \sum_{(k,\overline{\lambda}) \atop k \geq 1,(k) \cup \overline{\lambda}=\lambda}
\frac{1}{q^{k-1+2 \sum_{i=1}^{k-1} i \overline{m_i} + (2k-1) \sum_{i=k}^{n-k}
\overline{m_i}} q^{2 \sum_{i<j} i \overline{m_i} \overline{m_j} + \sum_i i
\overline{m_i}^2} \prod_i (\frac{1}{q})_{\overline{m_i}}}
\end{eqnarray*}

	Clearly if $(k) \cup \overline{\lambda}=\lambda$ then
$m_i=\overline{m_i}$ for $i \neq k$ and $m_k=\overline{m_k}+1$. Elementary
combinatorics then shows that this last expression is equal to

\begin{eqnarray*}
& & \sum_{(k,\overline{\lambda}) \atop k \geq 1,(k) \cup \overline{\lambda}=\lambda}
\frac{1}{q^{2 \sum_{i<j} i m_i m_j + \sum_i i m_i^2 - \sum_{i \geq k}
m_i} \prod_i (\frac{1}{q})_{\overline{m_i}}}\\
& = & \frac{1}{q^{2 \sum_{i<j} i m_i m_j + \sum_i i m_i^2} \prod_i
(\frac{1}{q})_{m_i}} \sum_{(k,\overline{\lambda}) \atop k \geq 1,(k) \cup \overline{\lambda}=\lambda} (1-1/q^{m_k}) q^{\lambda_k'}\\
& = & \frac{1}{q^{2 \sum_{i<j} i m_i m_j + \sum_i i m_i^2} \prod_i
(\frac{1}{q})_{m_i}} \sum_{k \geq 1} (1-1/q^{m_k}) q^{\lambda_k'}\\
& = & \frac{1}{q^{2 \sum_{i<j} i m_i m_j + \sum_i i m_i^2} \prod_i
(\frac{1}{q})_{m_i}} \sum_{k \geq 1} (q^{\lambda_k'}-q^{\lambda_{k+1}'})\\
& = &  (q^{\lambda_1'}-1) \frac{1}{q^{2 \sum_{i<j} i m_i m_j + \sum_i i m_i^2} \prod_i
(\frac{1}{q})_{m_i}} \\
& = & (q^{\lambda_1'}-1) \frac{1}{\prod_i q^{(\lambda_i')^2} (\frac{1}{q})_{m_i}},
\end{eqnarray*} where the final equality is Lemma \ref{useful}. This completes the proof of part 1 of the theorem. Part 2 follows from part 1 and Theorem \ref{classes}. \end{proof}

	Given Theorem \ref{probabform} it is now straightforward to
write down the corresponding cycle index generating functions. These
will be applied in Section \ref{S5}.

\begin{cor} \label{cycleindex} Let $x_{\phi,\lambda}$ be a collection
of variables. Let $\{\lambda_{\phi}(\alpha)\}$ be the rational
canonical form data of any $\alpha \in GL(n+1,q)$. Then

\begin{enumerate}

\item
\begin{eqnarray*}
&& \sum_{n=0}^{\infty} \frac{u^n}{|A(n,q)|} \sum_{\alpha \in
A(n,q)} \prod_{\phi \neq z} x_{\phi,\lambda_{\phi}(\alpha)}\\
& = &
(\sum_{\lambda: |\lambda|>0} \frac{x_{z-1,\lambda} u^{|\lambda|-1}
(q^{\lambda_1'}-1)}{\prod_i q^{(\lambda_i')^2}
(\frac{1}{q})_{m_i(\lambda)}}) \prod_{\phi \neq z,z-1} (\sum_{\lambda}
\frac{x_{\phi,\lambda} u^{|\lambda| \cdot deg(\phi)}}{\prod_i q^{deg(\phi) \cdot
(\lambda_i')^2} (\frac{1}{q^{deg(\phi)}})_{m_i(\lambda)}}).
\end{eqnarray*}

\item
\begin{eqnarray*}
&& \sum_{n=0}^{\infty} \frac{u^n}{|P(n,q)|} \sum_{\alpha \in
P(n,q)} \prod_{\phi \neq z} x_{\phi,\lambda_{\phi}(\alpha)}\\
& = & \sum_{a
\in F_q^*} (\sum_{\lambda: |\lambda|>0} \frac{x_{z-a,\lambda}
u^{|\lambda|-1} (q^{\lambda_1'}-1)}{(q-1) \prod_i q^{(\lambda_i')^2}
(\frac{1}{q})_{m_i(\lambda)}}) \prod_{\phi \neq z,z-a} (\sum_{\lambda}
\frac{x_{\phi,\lambda} u^{|\lambda| \cdot deg(\phi)}}{\prod_i q^{deg(\phi) \cdot
(\lambda_i')^2} (\frac{1}{q^{deg(\phi)}})_{m_i(\lambda)}}).
\end{eqnarray*}

\end{enumerate}

\end{cor}

\section{A Measure on Partitions and Probabilistic Algorithms}
\label{S3}

	To begin we define a probability measure which will be the
object of study in this section. Recall that $n(\lambda)$ is the
quantity $\sum_{i \geq 1} (i-1) \lambda_i$.

{\bf Definition:} For $0<u<1$ and $q>1$ the measure $N_{u,q}$ on the
set of all partitions of all positive integers is defined by

\begin{eqnarray*}
N_{u,q}(\lambda) & = & \prod_{r=1}^{\infty}(1-\frac{u}{q^r})
\frac{u^{|\lambda|-1} (q^{\lambda_1'}-1) }{\prod_i q^{(\lambda_i')^2}
(\frac{1}{q})_{m_i(\lambda)}}\\
& = & \prod_{r=1}^{\infty}(1-\frac{u}{q^r}) \frac{u^{|\lambda|-1}
(q^{\lambda_1'}-1)
P_{\lambda}(\frac{1}{q},\frac{1}{q^2},\cdots;\frac{1}{q})}{q^{n(\lambda)}}.
\end{eqnarray*} The equivalence between the two definitions of
$N_{u,q}(\lambda)$ follows from Lemma \ref{useful}.

	Before proving that $N_{u,q}$ is indeed a probability measure,
we recall the analogous measure $M_{u,q}$ on the set of all partitions
of all natural numbers. The survey \cite{F3} summarizes what is known
about $M_{u,q}$.

\begin{eqnarray*}
M_{u,q}(\lambda) & = & \prod_{r=1}^{\infty}(1-\frac{u}{q^r})
\frac{u^{|\lambda|} }{\prod_i q^{(\lambda_i')^2}
(\frac{1}{q})_{m_i(\lambda)}}\\
& = & \prod_{r=1}^{\infty}(1-\frac{u}{q^r}) \frac{u^{|\lambda|}
P_{\lambda}(\frac{1}{q},\frac{1}{q^2},\cdots;\frac{1}{q})}{q^{n(\lambda)}}.
\end{eqnarray*}

	Lemma \ref{ismeasure} proves that $N_{u,q}$ is a probability
measure. We include two proofs--one using an identity of Macdonald
from symmetric function theory and another using the cycle index of
$A(n,q)$ and the fact that $M_{u,q}$ is a probability measure for
$q>1$ and $0<u<1$. (Two other proof methods are to use either
Steinberg's result that there are $q^{n^2}$ unipotent elements in
$A(n,q)$ together with the cycle index or else to use the Markov chain
construction later in this section together with an identity of
Cauchy).

\begin{lemma} \label{ismeasure} If $q>1$ and $0<u<1$ then $N_{u,q}$
defines a probability measure. \end{lemma}

\begin{proof} (First proof) Equation 3 on page 219 of \cite{Mac}
states that

\[ \prod_{i \geq 1} (1+x_iy)/(1-x_i) = \sum_{\lambda} t^{n(\lambda)}
\prod_{j=1}^{\lambda_1'} (1+t^{1-j}y) P_{\lambda}(x;t) .\] The result
now follows by first setting $x_i=u/q^i,t=1/q$ and taking coefficients of
$y$ on both sides, and then using the fact that
$P_{\lambda}(u/q,u/q^2,\cdots;1/q)=u^{|\lambda|}P_{\lambda}(1/q,1/q^2,\cdots;1/q)$.
\end{proof}

\begin{proof} (Second proof) Setting all $x_{\phi,\lambda}$ in the
cycle index of $A(n,q)$ equal to 1 gives the identity

\[ \frac{1}{1-u} = (\prod_{r=1}^{\infty} \frac{1}{1-u/q^r}
\sum_{\lambda: |\lambda|>0} N_{u,q}(\lambda)) \prod_{\phi \neq z,z-1}
(\prod_{r=1}^{\infty} \frac{1}{1-u^{deg(\phi)}/q^{r \cdot deg(\phi)}})
\sum_{\lambda} M_{u^{deg(\phi)},q^{deg(\phi)}}(\lambda).\] Since
$M_{u,q}$ is a probability measure it follows that

\[ \frac{1}{1-u} = (\prod_{\phi \neq z} \prod_{r=1}^{\infty}
\frac{1}{1-u^{deg(\phi)}/q^{r \cdot deg(\phi)}}) \sum_{\lambda: |\lambda|>0}
N_{u,q}(\lambda).\] Setting all variables in the cycle index for
$GL(n,q)$ equal to 1 and using the fact that $M_{u,q}$ is a
probability measure implies that

\[ \frac{1}{1-u} = \prod_{\phi \neq z} \prod_{r=1}^{\infty}
\frac{1}{1-u^{deg(\phi)}/q^{r \cdot deg(\phi)}}\] (of course this can be proved
directly). The lemma follows. \end{proof}

	We require an elementary lemma about Taylor series. For its
statement $[u^n]f(u)$ denotes the coefficient of $u^n$ in a polynomial
$f(u)$.

\begin{lemma} \label{largen} If $f(1)<\infty$ and the Taylor series of
$f$ around 0 converges at $u=1$, then
\[ lim_{n \rightarrow \infty} [u^n] \frac{f(u)}{1-u} = f(1). \]
\end{lemma}

\begin{proof} Write the Taylor expansion $f(u) = \sum_{n=0}^{\infty}
a_n u^n$. Then observe that $[u^n] \frac{f(u)}{1-u} = \sum_{i=0}^n
a_i$.  \end{proof}

	Theorem \ref{interpret} shows that the measure $N_{u,q}$ is a
fundamental object for understanding the probability theory of
$A(n,q)$. We remark that the idea behind it--auxiliary
randomization--is a mainstay of statistical mechanics.

\begin{theorem} \label{interpret}

\begin{enumerate}

\item Fix $u$ with $0<u<1$. Then choose a random number $N$ with
probability of getting $n$ equal to $(1-u)u^n$. Choose $\alpha$
uniformly in $A(N,q)$. Then as $\phi$ varies any finite number of the random partitions
$\lambda_{\phi}(\alpha)$ are independent random variables, with
$\lambda_{z-1}$ distributed according to the measure $N_{u,q}$ and all
other $\lambda_{\phi}$ distributed according to the measure
$M_{u^{deg(\phi)},q^{deg(\phi)}}$.

\item Choose $\alpha$ uniformly in $A(n,q)$. Then as $n \rightarrow
\infty$, any finite number of the random partitions $\lambda_{\phi}(\alpha)$ are
independent random variables, with $\lambda_{z-1}$ distributed
according to the measure $N_{1,q}$ and all other $\lambda_{\phi}$
distributed according to the measure $M_{1,q^{deg(\phi)}}$.

\item Fix $u$ with $0<u<1$. Then choose a random number $N$ with
probability of getting $n$ equal to $(1-u)u^n$. Choose $\alpha$
uniformly in $P(N,q)$. Then as $\phi$ varies over polynomials of
degree at least two, any finite number of the random partitions $\lambda_{\phi}(\alpha)$
are independent random variables distributed according to the measure
$M_{u^{deg(\phi)},q^{deg(\phi)}}$. The partitions $\lambda_{z-a}$ are
independent of $\lambda_{\phi}$ for $deg(\phi) \geq 2$ but depend on
each other; any particular $\lambda_{z-a}$ has as its distribution the
mixture $\frac{N_{u,q}}{q-1} + \frac{(q-2) M_{u,q}}{q-1}$.

\item Choose $\alpha$ uniformly in $P(n,q)$. Then as $n \rightarrow
\infty$ and $\phi$ varies over polynomials of degree at least two, any finite number of the
random partitions $\lambda_{\phi}(\alpha)$ are independent random
variables distributed according to the measure
$M_{1,q^{deg(\phi)}}$. The partitions $\lambda_{z-a}$ are independent
of $\lambda_{\phi}$ for $deg(\phi) \geq 2$ but depend on each other;
any particular $\lambda_{z-a}$ has as its distribution the mixture
$\frac{N_{1,q}}{q-1} + \frac{(q-2) M_{1,q}}{q-1}$.

\end{enumerate}
\end{theorem}

\begin{proof} As explained in the second proof of Lemma
\ref{ismeasure}, we know that \[ \frac{1}{1-u} = \prod_{\phi \neq z}
\prod_{r=1}^{\infty} \frac{1}{1-u^{deg(\phi)}/q^{r \cdot deg(\phi)}}.\] Multiplying
this by the cycle index of $A(n,q)$ gives that

\[ \sum_{n=0}^{\infty} \frac{(1-u) u^n}{|A(n,q)|} \sum_{\alpha \in
A(n,q)} \prod_{\phi \neq z} x_{\phi,\lambda_{\phi}(\alpha)} =
(\sum_{\lambda: |\lambda|>0} x_{z-1,\lambda} N_{u,q}(\lambda))
\prod_{\phi \neq z,z-1} (\sum_{\lambda} x_{\phi,\lambda}
M_{u^{deg(\phi)},q^{deg(\phi)}}(\lambda)).\] This proves the first
assertion of the theorem. For the second assertion use Lemma
\ref{largen}. The proofs of the third and fourth assertions are
almost identical so are omitted.
\end{proof}

	The remainder of this section considers probabilistic methods
for growing random partitions according to the measure $N_{u,q}$,
analogous to those for $M_{u,q}$ in \cite{F0},\cite{F1}. For this
recall that a standard Young tableau $T$ of size $n$ is a partition of
$n$ with each box containing one of $\{1,\cdots,n\}$ such that each
of $\{1,\cdots,n\}$ appears exactly once and the numbers increase in
each row and column of $T$. For instance,

\[ \begin{array}{c c c c c}
                \framebox{1} & \framebox{3} & \framebox{5} & \framebox{6} &   \\
                \framebox{2} & \framebox{4} & \framebox{7} & &    \\
                \framebox{8} & \framebox{9} &  &  &    
          \end{array} \] is a standard Young tableau.

	Our first method for growing partitions according to the
measure $N_{u,q}$ is an algorithm we call the Affine Young Tableau
Algorithm, so as to distinguish it from the Young Tableau Algorithm
for growing partitions according to the measure $M_{u,q}$. For its
statement, we need the Young Tableau Algorithm \cite{F0},\cite{F1}.

\begin{center}
The Young Tableau Algorithm
\end{center}

\begin{description}

\item [Step 0] Start with $N=1$ and $\lambda$ the empty
partition. Also start with a collection of coins indexed by the
natural numbers, such that coin $i$ has probability $\frac{u}{q^i}$ of
heads and probability $1-\frac{u}{q^i}$ of tails.

\item [Step 1] Flip coin $N$.

\item [Step 2a] If coin $N$ comes up tails, leave $\lambda$ unchanged,
set $N=N+1$ and go to Step 1.

\item [Step 2b] If coin $N$ comes up heads, choose an integer $S>0$
according to the following rule. Set $S=1$ with probability $\frac
{q^{N-\lambda_1'}-1} {q^N-1}$. Set $S=s>1$ with probability
$\frac{q^{N-\lambda_s'}-q^{N-\lambda_{s-1}'}}{q^N-1}$. Then increase
the size of column $s$ of $\lambda$ by 1 and go to Step 1.

\end{description}

        As an example of the Young Tableau Algorithm, suppose we are
at Step 1 with $\lambda$ equal to the following partition:
        
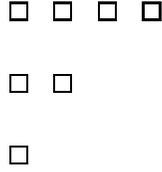
\[ \begin{array}{c c c c}
                \framebox{} & \framebox{} & \framebox{} & \framebox{}  \\
                \framebox{} & \framebox{} &  &      \\
                \framebox{} &  &  &    \\
                 & & &  
          \end{array} \] Suppose also that $N=4$ and that coin 4 had
already come up heads once, at which time we added to column 1, giving
$\lambda$. Now we flip coin 4 again and get heads, going to Step
2b. We add to column $1$ with probability $\frac{q-1}{q^4-1}$, to
column $2$ with probability $\frac{q^2-q}{q^4-1}$, to column $3$ with
probability $\frac{q^3-q^2} {q^4-1}$, to column $4$ with probability
$0$, and to column $5$ with probability $\frac{q^4-q^3}{q^4-1}$. We
then return to Step 1.

\begin{theorem} (\cite{F0},\cite{F1}) \label{TableauAlg} For $0<u<1$
and $q>1$, the Young Tableau Algorithm generates partitions which are
distributed according to the measure $M_{u,q}$. \end{theorem}

	Next we point out that the Young Tableau Algorithm can easily
be made to terminate. This idea was developed in joint work with Mark
Huber surveyed in \cite{F3} and goes as follows. Let $a_N$ be the
number of times that coin $N$ comes up heads; the idea is simply to
first determine the random vector $(a_1,a_2,\cdots)$ and then grow the
partitions as in Step 2b of the Young Tableau Algorithm. So let us
explain how to determine $(a_1,a_2,\cdots)$. For $N \geq 1$ let
$t^{(N)}$ be the probability that all tosses of all coins numbered $N$
or greater are tails. For $N \geq 1$ and $j \geq 0$ let $t^{(N)}_j$ be
the probability that some toss of a coin numbered $N$ or greater is a
head and that coin $N$ comes up heads $j$ times. It is simple to write
down expressions for $t^{(N)},t^{(N)}_0,t^{(N)}_1,\cdots$ and clearly
$t^{(N)}+\sum_{j \geq 0} t^{(N)}_j=1$. The basic operation a computer
can perform is to produce a random variable $U$ distributed uniformly
in the interval $[0,1]$. By dividing $[0,1]$ into intervals of length
$t^{(1)},t^{(1)}_0,t^{(1)}_1,\cdots$ and seeing where $U$ is located,
one arrives at the value of $a_1$. Furthermore, if $U$ landed in the
interval of length $t^{(1)}$ then all coins come up tails and the
vector $(a_1,a_2,\cdots)$ is determined. Otherwise, move on to coin 2,
dividing $[0,1]$ into intervals of length
$t^{(2)},t^{(2)}_0,t^{(2)}_1,\cdots$ and so on.

	Now we describe the Affine Young Tableau Algorithm.

\begin{center}
The Affine Young Tableau Algorithm
\end{center}

\begin{description}

\item [Step 1] Run the Young Tableau Algorithm so as to generate a
partition $\lambda$ distributed as $M_{u,q}$.

\item [Step 2] Set $S=1$ with probability $\frac{1}{q^{\lambda_1'}}$
and $S=s>1$ with probability $\frac{1}{q^{\lambda_s'}} -
\frac{1}{q^{\lambda_{s-1}'}}$. Then increase the size of column $s$ of
$\lambda$ by $1$.

\end{description}

	Theorem \ref{affineOK} shows that the Affine Young Tableau
Algorithm grows partitions distributed as $N_{u,q}$. For its proof, we
use a different method for generating $M_{u,q}$. Recall that the Young
lattice is the set of all partitions of all natural numbers, with a
directed edge drawn from the partition $\lambda$ to partition
$\Lambda$ if the diagram of $\lambda$ is contained in the diagram of
$\Lambda$ and $|\Lambda|=|\lambda|+1$.

\begin{theorem} \label{weights} (\cite{F1}) Put weights
$m_{\lambda,\Lambda}$ on the edges of Young lattice according to the
rules:

\begin{enumerate}

\item $m_{\lambda,\Lambda} =
  \frac{u}{q^{\lambda_1'} (q^{\lambda_1'+1}-1)}$ if the diagram of
$\Lambda$ is obtained from that of $\lambda$ by adding a box to column
1.

\item $m_{\lambda,\Lambda} = \frac{u(q^{-\lambda_s'}-q^{-
\lambda_{s-1}'})}{q^{\lambda_1'}-1}$ if the diagram of $\Lambda$ is
obtained from that of $\lambda$ by adding a box to column $s>1$.

\end{enumerate} Then the following formula holds:

\[ M_{u,q}(\lambda) = \left[\prod_{r=1}^{\infty}
(1-\frac{u}{q^r})\right] \sum_{\gamma} \prod_{i=0}^{|\lambda|-1}
m_{\gamma_i,\gamma_{i+1}} \] where the sum is over all directed paths
$\gamma$ from the empty partition to $\lambda$, and the $\gamma_i$ are
the partitions along the path $\gamma$. Furthermore, if the Young
tableau $T$ corresponds to the path $\gamma$ then the probability that
the Young Tableau Algorithm outputs $T$ is equal to \[
\left[\prod_{r=1}^{\infty} (1-\frac{u}{q^r})\right]
\prod_{i=0}^{|\lambda|-1} m_{\gamma_i,\gamma_{i+1}}.\] \end{theorem}

	Now we can show that the Affine Young Tableau Algorithm works.
 
\begin{theorem} \label{affineOK} The Affine Young Tableau Algorithm
grows partitions distributed as $N_{u,q}$. \end{theorem}

\begin{proof} Let $\lambda^{(s)}$ denote the the shape obtained by
decreasing the size of column $s$ of $\lambda$ by one and let
$M_{u,q}(\lambda^{(s)})=0$ if $\lambda^{(s)}$ is not a partition. To
prove the theorem it is sufficient to show that for $|\lambda|>0$,

\[ \frac{1}{q^{\lambda_1'-1}} M_{u,q}(\lambda^{(1)}) + \sum_{s>1}
(\frac{1}{q^{\lambda_s'-1}} - \frac{1}{q^{\lambda_{s-1}'}})
M_{u,q}(\lambda^{(s)}) = N_{u,q}(\lambda) .\] Theorem \ref{weights} shows
that 

\begin{eqnarray*}
& &  \frac{1}{q^{\lambda_1'-1}} M_{u,q}(\lambda^{(1)}) + \sum_{s>1}
(\frac{1}{q^{\lambda_s'-1}} - \frac{1}{q^{\lambda_{s-1}'}})
M_{u,q}(\lambda^{(s)})\\
& = & \prod_{r=1}^{\infty} (1-\frac{u}{q^r}) [ \sum_{s \geq 1} \sum_{\gamma: \emptyset \rightarrow \lambda^{(s)} \rightarrow \lambda } (\prod_{i=0}^{|\lambda|-2} m_{\gamma_i,\gamma_{i+1}}) (\frac{(q^{\lambda_1'}-1) m_{\lambda^{(s)},\lambda}}{u})]
\end{eqnarray*} where the inner sum is over all paths in the Young lattice from the empty partition to $\lambda^{(s)}$ to $\lambda$ and all other notation is as in Theorem \ref{weights}. Simplifying further and using Theorem \ref{weights} again gives that

\[ \prod_{r=1}^{\infty} (1-\frac{u}{q^r}) [\frac{q^{\lambda_1'}-1}{u}
\sum_{\gamma: \emptyset \rightarrow \lambda} \prod_{i=0}^{|\lambda|-1}
m_{\gamma_i,\gamma_{i+1}}] = \frac{q^{\lambda_1'}-1}{u} M_{u,q}(\lambda)
= N_{u,q}(\lambda).\]

\end{proof}

{\bf Remarks:}
\begin{enumerate}

\item Theorem \ref{affineOK} gives another proof that $N_{u,q}$ is a
probability measure.

\item From Theorem \ref{weights} one has an explicit formula for the
chance that the Young Tableau Algorithm outputs a given Young tableau
(page 565 of \cite{F1}). Since the Affine Young Tableau Algorithm
randomly adds on one additional box, one can write down a formula for
the chance that it generates a given Young Tableau.

\item Two other probabilistic algorithms are given for growing
partitions according to $M_{u,q}$ on pages 581 and 585 of
\cite{F1}. These yield algorithms for sampling from $N_{u,q}$ simply
by adding an additional box as in Step 2 of the Affine Young Tableau
Algorithm.

\end{enumerate}

	Next we describe a rather surprising method for sampling from
$N_{u,q}$ using Markov chains. We use the notation that
$(\frac{u}{q})_i$ is equal to $(1-u/q) \cdots (1-u/q^i)$. $Prob.(E)$
will denote the probability of an event $E$ under the measure
$N_{u,q}$.

\begin{theorem} \label{markovgl} Starting with $\lambda_1'=a \geq 1$
with probability $Q(a)=\frac{(\prod_{r=1}^{\infty} (1-u/q^r))
u^{a-1}}{q^{a^2-a} (\frac{u}{q})_a (\frac{1}{q})_{a-1}}$, define in
succession $\lambda_2',\lambda_3',\cdots$ according to the rule that
if $\lambda_i'=a$, then $\lambda_{i+1}'=b$ with probability \[ K(a,b)
= \frac{u^b (\frac{1}{q})_a (\frac{u}{q})_a}{q^{b^2}
(\frac{1}{q})_{a-b} (\frac{1}{q})_b (\frac{u}{q})_b}.\] Then the
resulting partition is distributed according to
$N_{u,q}$. \end{theorem}

\begin{proof} The $N_{u,q}$ probability of choosing a partition
with $\lambda_i'=r_i$ for all $i$ is

\[ Prob.(\lambda_1'=r_1) \prod_{i=1}^{\infty} \frac{Prob.(
\lambda_1'=r_1,\cdots,\lambda_{i+1}'=r_{i+1})} {Prob.(\lambda_1'=
r_1,\cdots,\lambda_{i}'=r_{i})}.\] Theorem \ref{fixedspace} of Section
\ref{S4} shows that $Prob.(\lambda_1'=r_1)=Q(r_1)$. Thus it is enough
to prove that

\[ \frac{Prob.( \lambda_1'=r_1,\cdots,\lambda_{i-1}'=r_{i-1},
\lambda_i'=a, \lambda_{i+1}'=b)} {Prob.(\lambda_1'=
r_1,\cdots,\lambda_{i-1}'=r_{i-1}, \lambda_{i}'=a)} = \frac{u^b
(\frac{1}{q})_a (\frac{u}{q})_a}{q^{b^2} (\frac{1}{q})_{a-b}
(\frac{1}{q})_b (\frac{u}{q})_b},\] for all $i \geq
1,a,b,r_1,\cdots,r_{i-1}$.
 
	Letting $T(a)$ be the $M_{u,q}$ probability that
$\lambda_1'=a$ it is proved in \cite{F1} that $T(a)=\frac{u^a
(\frac{u}{q})_{\infty}}{q^{a^2} (\frac{1}{q})_a
(\frac{u}{q})_a}$. Next one calculates that for $i>2$

\[ \sum_{\lambda: \lambda_1'=r_1,\cdots,\lambda_{i-1}'=r_{i-1} \atop
\lambda_i'=a} N_{u,q} (\lambda) = \frac{(q^{r_1}-1)
u^{r_1+\cdots+r_{i-1}-1}} {q^{r_1^2+\cdots+r_{i-1}^2}
(\frac{1}{q})_{r_1-r_2} \cdots (\frac{1}{q})_{r_{i-2}-r_{i-1}}
(\frac{1}{q})_{r_{i-1}-a}} T(a).\] Similarly, observe that \[
\sum_{\lambda: \lambda_1'=r_1,\cdots,\lambda_{i-1}'=r_{i-1} \atop
\lambda_i'=a,\lambda_{i+1}'=b} N_{u,q}(\lambda) = \frac{(q^{r_1}-1)
u^{r_1+\cdots+r_{i-1}+a-1}} {q^{r_1^2+\cdots+r_{i-1}^2+a^2}
(\frac{1}{q})_{r_1-r_2} \cdots (\frac{1}{q})_{r_{i-2}-r_{i-1}}
(\frac{1}{q})_{r_{i-1}-a} (\frac{1}{q})_{a-b}} T(b).\] Thus the ratio
of these two expressions is \[ \frac{u^b (\frac{1}{q})_a
(\frac{u}{q})_a}{q^{b^2} (\frac{1}{q})_{a-b} (\frac{1}{q})_b
(\frac{u}{q})_b}, \] as desired. The case $i=1$ must be checked
separately but the same ratio results.
\end{proof}

	The algorithm of Theorem \ref{markovgl} runs on a computer
because of the well known fact that to sample from a discrete
distribution one divides $[0,1]$ into intervals of the appropriate
lengths and generates a uniform variable in $[0,1]$.

	To conclude the section we explain how to sample from
$N_{u,q}$ given that the size of the partition is $n+1$, which is
clearly the same as sampling the partition $\lambda_{z-1}$ for a
unipotent element of $A(u,q)$. One naive method is to keep picking
from $N_{u,q}$ until one obtains a partition of size $n$; however this
is very slow and not theoretically useful.

\begin{center}
Algorithm for Sampling from $N_{u,q}$ given that $|\lambda|=n+1$
\end{center}

\begin{description}

\item [Step 0] Start with $N=1$ and $\lambda$ the empty partition.

\item [Step 1] If $n=0$ then go to step 3. Otherwise set $h=1-\frac{1}{q^n}$.

\item [Step 2] Flip a coin with probability of heads $h$.

\item [Step 2a] If the toss of Step 2 came up tails, increase the value
of $N$ by $1$ and go to Step 2.

\item [Step 2b] If the toss of Step 2 comes up heads, decrease the value
of $n$ by $1$, increase $\lambda$ according to the rule of Step 2b of
the Young Tableau Algorithm (which depends on $N$), and then go to Step 1.

\item [Step 3] Perform Step 2 of the Affine Young Tableau Algorithm.

\end{description}

	Before proving Theorem \ref{unipotentsample}, we remark that
it is an adaptian of an analogous result for $GL(n,q)$ which was
obtained in joint work with Mark Huber, surveyed in \cite{F3}.

\begin{theorem} \label{unipotentsample} The algorithm for sampling
from $N_{u,q}$ given that $|\lambda|=n+1$ is valid. \end{theorem}

\begin{proof} The survey \cite{F3} proved that if one replaced Step 3
in the above algorithm by stopping then one would sample from
$M_{u,q}$ given that $|\lambda|=n$. Let $\lambda^{(s)}$ denote
$\lambda$ with the size of column $s$ decreased by one and let $E(s)$
be the probability that Step 2 of the Affine Young Tableau Algorithm
adds to column $s$ of $\lambda^{(s)}$. The result follows because

\begin{eqnarray*}
\frac{N_{u,q}(\lambda)}{\sum_{\lambda: |\lambda|=n+1} N_{u,q}(\lambda)} & = & 
\frac{N_{u,q}(\lambda)}{\sum_{\lambda: |\lambda|=n} M_{u,q}(\lambda)}\\
& = & \frac{\sum_s M_{u,q}(\lambda^{(s)}) E(s)}{\sum_{\lambda: |\lambda|=n} M_{u,q}(\lambda)}\\
& = & \sum_s E(s) \frac{M_{u,q}(\lambda^{(s)})}{\sum_{\lambda: |\lambda|=n} M_{u,q}(\lambda)}.
\end{eqnarray*}

\end{proof}

\section{Applications of Probabilistic Algorithms} \label{S4}

	This section consider applications of the probabilistic
algorithms of Section \ref{S3}. The results here parallel those of
\cite{F1} for $GL(n,q)$ but are genuinely different. We do omit the
formula for the chance that an element of $A(n,q)$ has a given
characteristic polynomial, as this follows from the $GL(n,q)$ result
and offers no new insights. Furthermore we only state the results for
$A(n,q)$ as the extensions to $P(n,q)$ are trivial. This section uses
the notation that $[u^n] f(u)$ denotes the coefficient of $u^n$ in a
polynomial $f(u)$.

	Let $P_{A,n}(k,q)$ be the chance that an element of $A(n,q)$
has a $k$ dimensional fixed space and let $P_{A,\infty}(k,q)$ be the
$n \rightarrow \infty$ limit of $P_{A,n}(k,q)$. Theorem
\ref{fixedspace} will show that

\[ P_{A,n}(k,q) = \frac{1}{|A(k-1,q)|} \sum_{i=0}^{n-k+1}
\frac{(-1)^i}{q^{ki} (q^i-1) \cdots (q-1)} \] and

\[ P_{A,\infty}(k,q) = [ \prod_{r=1}^{\infty} (1-\frac{1}{q^r}) ]
\frac{(1/q)^{k^2-k}}{(1-1/q)^2 \cdots (1-1/q^{k-1})^2 (1-1/q^k)},\]
where we use the notation that $|A(-1,q)|=0$.  Furthermore, a
probabilistic interpretation will be given to the products in
$P_{A,\infty}(k,q)$. The first step is to connect the chance that an
element $\alpha$ of $A(n,q)$ has a $k$ dimensional fixed space with
the partitions in the rational canonical form of $\alpha$.

\begin{lemma} \label{RatParts} (\cite{F1}) The dimension of the fixed
space of an element $\alpha$ of $GL(n,q)$ is equal to
$\lambda_{z-1,1}'(\alpha)$ (i.e. the number of parts of the partition
corresponding to the polynomial $z-1$ in the rational canonical form
of $\alpha$).  \end{lemma}

	To proceed further, some notation is necessary. Let $T$ be a
standard Young tableau with $k$ parts. Let $T_{(i,j)}$ be the entry in
row $i$ and column $j$ of $T$. We define numbers
$h_1(T),\cdots,h_k(T)$ associated with $T$. Let
$h_m(T)=T_{(m+1,1)}-T_{(m,1)}-1$ for $1 \leq m \leq k-1$ and let
$h_k(T)=|T|-T_{(k,1)}$. So if $k=3$ and $T$ is the tableau

\[ \begin{array}{c c c c c}
                \framebox{1} & \framebox{3} & \framebox{5} & \framebox{6} &   \\
                \framebox{2} & \framebox{4} & \framebox{7} &  &    \\
                \framebox{8} & \framebox{9} &  &  &    
          \end{array} \] then $h_1(T)=2-1-1=0$, $h_2(T)=8-2-1=5$, and
$h_3(T)=9-8=1$. View $T$ as being created by the Affine Young Tableau
Algorithm. Then for $1 \leq m \leq k-1$, $h_m(T)$ is the number of
boxes added to $T$ after it becomes a tableau with $m$ parts and before
it becomes a tableau with $m+1$ parts. $h_k(T)$ is the number of boxes
added to $T$ after it becomes a tableau with $k$ parts. The proof of
Theorem $\ref{Interp}$ will show that if one conditions $T$ chosen
from the measure $N_{u,q}$ on having $k$ parts, then the random
variables $h_1(T),\cdots,h_{k}(T)$ are independent, where
$h_1(T),\cdots,h_{k}(T)$ are geometric with parameters $\frac{u}{q},
\cdots,\frac{u}{q^{k}}$. 

\begin{theorem} \label{Interp}

\[ \sum_{\lambda: \lambda_1'=k} N_{u,q}(\lambda) = \frac{u^{k-1} }{|A(k-1,q)|} \frac{\prod_{r=1}^{\infty}
(1-\frac{u}{q^r})}{\prod_{r=1}^k (1-\frac{u}{q^r})} \]

\end{theorem}

\begin{proof} We sum over all Young tableaux $T$ with $k$ parts the chance that the Affine Young Tableau
Algorithm outputs $T$. Recall from the second remark after Theorem
\ref{weights} that one can compute the probability that the Affine
Young Tableau Algorithm outputs any given Young tableau. The point is
that one can in fact compute the probability that the Affine Young
Tableau Algorithm produces a tableau $T$ with given values of the
$h$'s.

        To do this, consider what happens during Step 1 of the Affine Young
Tableau Algorithm. Suppose that one moves along the Young lattice from
a partition with $m$ parts. Theorem $\ref{weights}$ implies that the
weight for adding to column $1$ is $\frac{u}{q^m(q^{m+1}-1)}$, and
that the sum of the weights for adding to some other column is
$\frac{u}{q^m}$. Thus the chance that Step 1 of the Affine Young
Tableau Algorithm yields a tableau with given values
$h_1,\cdots,h_{k-1}$ and all other $h_i$=0 is

\[ \prod_{r=1}^{\infty} (1-\frac{u}{q^r}) \frac{u^{k-1}}{|GL(k-1,q)|}
\prod_{m=1}^{k-1} (\frac {u}{q^m})^{h_m}. \] 

	In order for the Affine Young Tableau Algorithm to generate a
partition with $m$ parts, either Step 1 generates a partition with
$m-1$ parts and Step 2 increases the size of column $1$ or else Step 1
generates a partition with $m$ parts and Step 2 increases the size of
column $s>1$. Thus the probability that the Affine Young Tableau
Algorithm generates a partition with $k$ parts is

\begin{eqnarray*}
& & \sum_{\lambda: \sum \lambda_1'=k} N_{u,q}(\lambda)\\
& = & 
\left[\prod_{r=1}^{\infty} (1-\frac{u}{q^r})\right] \frac{u^{k-1}}{|GL(k-1,q)|}
\prod_{m=1}^{k-1} \sum_{h_m=0}^{\infty} (\frac{u}{q^m})^{h_m} [\frac{1}{q^{k-1}} + \frac{u}{q^{k-1}(q^k-1)} \sum_{h_k \geq 0} (\frac{u}{q^k})^{h_k} (1-\frac{1}{q^k})]\\
& = & \left[\prod_{r=1}^{\infty} (1-\frac{u}{q^r})\right] \frac{u^{k-1}}{q^{ k-1} |GL(k-1,q)|}
\prod_{m=1}^{k-1} \sum_{h_m=0}^{\infty} (\frac{u}{q^m})^{h_m} [\sum_{h_k \geq 0} (\frac{u}{q^k})^{h_k}]\\
& = & \frac{u^{k-1} }{|A(k-1,q)|} \frac{\prod_{r=1}^{\infty}
(1-\frac{u}{q^r})}{\prod_{r=1}^k (1-\frac{u}{q^r})}.
\end{eqnarray*}
\end{proof}

	One further lemma will be used.

\begin{lemma} \label{Gold} (\cite{GR})

\[ \prod_{r=1}^{\infty} (1-\frac{u}{q^r}) = \sum_{i=0}^{\infty}
\frac{(-u)^i}{(q^i-1) \cdots (q-1)} \]

\end{lemma}

	Now the formulas for $P_{A,n}(k,q)$ and $P_{A,\infty}(k,q)$
can be proved. Recall that we use the notation that $|A(-1,q)|=0$.

\begin{theorem} \label{fixedspace} 

\begin{enumerate}

\item \[ P_{A,n}(k,q) = \frac{1}{|A(k-1,q)|} \sum_{i=0}^{n-k+1}
\frac{(-1)^i}{q^{ki} (q^i-1) \cdots (q-1)}. \]

\item \[ P_{A,\infty}(k,q) = [ \prod_{r=1}^{\infty} (1-\frac{1}{q^r}) ]
\frac{(1/q)^{k^2-k}}{(1-1/q)^2 \cdots (1-1/q^{k-1})^2 (1-1/q^k)}.\]

\end{enumerate}
\end{theorem}

\begin{proof} In the equation of part one of Corollary \ref{cycleindex}
set $x_{z-1,\lambda}=1$ if $\lambda$ has k parts and
$x_{z-1,\lambda}=0$ otherwise. Also set $x_{\phi,\lambda}=1$ for $\phi \neq z-1$. It follows from Theorem
\ref{Interp} and Lemma \ref{Gold} that
\begin{eqnarray*}
P_{A,n}(k,q) & = & [u^n] \prod_{\phi \neq z} \prod_{r=1}^{\infty} (\frac{1}{1-\frac{u^{deg(\phi)}}{q^{r \cdot deg(\phi)}}}) \sum_{\lambda:\lambda_1'=k}
N_{u,q}(\lambda)\\
& = & [u^n] \frac{1}{1-u} \sum_{\lambda:\lambda_1'=k}
N_{u,q}(\lambda)\\
& = & [u^n] \frac{u^{k-1} \prod_{r=1}^{\infty} (1-\frac{u}{q^{k+r}})}{(1-u) |A(k-1,q)|}\\
& = & \frac{1}{|A(k-1,q)|} [u^{n-k+1}] \frac{1}{1-u}
\sum_{i=0}^{\infty} \frac{(-1)^i(uq^{-k})^i}{(q^i-1) \cdots (q-1)}\\
& = & \frac{1}{|A(k-1,q)|} \sum_{i=0}^{n-k+1}
\frac{(-1)^i}{q^{ki} (q^i-1) \cdots (q-1)}.
\end{eqnarray*}

	For the second assertion of the theorem use Lemma \ref{largen}
and Theorem \ref{Interp} to conclude that 
\begin{eqnarray*}
P_{A,\infty}(k,q) & = & lim_{n \rightarrow \infty}[u^n] \frac{1}{1-u} \sum_{\lambda:\lambda_1'=k} N_{u,q}(\lambda)\\
& = & \sum_{\lambda:\lambda_1'=k} N_{1,q}(\lambda)\\	
& = &  \frac{1}{|A(k-1,q)|} \frac{\prod_{r=1}^{\infty}
(1-\frac{1}{q^r})}{\prod_{r=1}^k (1-\frac{1}{q^r})}
\end{eqnarray*} as desired.
\end{proof}

	We remark that the analog of Theorem \ref{fixedspace} for
$GL(n,q)$ was first proved by Rudavlis and Shinoda \cite{RS}, using
Moebius inversion on the lattice of subspaces of a vector
space. Theorem \ref{fixedspace} can be proved along the same lines but
we omit the details as the argument is less insightful and incongruous
with the theme of this paper.

	The final result in this section is an enumeration of
unipotent elements in $A(n,q)$ of a given rank.

\begin{theorem} The number of unipotent elements of rank $k$ in
$A(n,q)$ is equal to

\[ \frac{|A(n,q)|}{|A(k-1,q)|} \frac{(1-1/q^k) \cdots
(1-1/q^n)}{q^{n-k+1} (1-1/q) \cdots (1-1/q^{n-k+1})}.\] \end{theorem}

\begin{proof} In the equation of part 2 of Corollary \ref{cycleindex}
set $x_{z-1,\lambda}=1$ if $\lambda$ has k parts and
$x_{\phi,\lambda}=0$ otherwise. Using Theorem \ref{Interp} one sees
that the sought number is

\begin{eqnarray*}
& & |A(n,q)| [u^n] \frac{u^{k-1}
}{|A(k-1,q)|}\frac{1}{\prod_{r=1}^k (1-\frac{u}{q^r})}\\
& = & \frac{|A(n,q)|}{|A(k-1,q)|} [u^{n-k+1}] \prod_{r=1}^k \frac{1}{1-\frac{u}{q^r}}\\
& = & \frac{|A(n,q)|}{|A(k-1,q)|} \frac{(1-1/q^k) \cdots
(1-1/q^n)}{q^{n-k+1} (1-1/q) \cdots (1-1/q^{n-k+1})}.
\end{eqnarray*} 
\end{proof}

\section{Applications of Cycle Indices} \label{S5}

	The purpose of this section is to apply the cycle index
generating functions of Section \ref{S2} to obtain precise estimates
for the probabilities that an element of $A(n,q)$ is separable,
cyclic, or semisimple (these terms will be defined as needed). As
these probabilities are identical for $P(n,q)$ results will only be
stated for $A(n,q)$. The section is organized by discussing separable,
cyclic, and semisimple probabilities, and in that order.

	The use of generating function methods to provide similar
estimates for the finite classical groups has appeared in earlier
papers (the paper \cite{F5} computes all three limits for $GL$,
\cite{W} computes the separable and cyclic limits and bounds the
convergence rates, and \cite{FNP} obtains convergence rates for the
semisimple case and extensions to the finite classical groups). The
reader unfamiliar with the cycle index of $GL(n,q)$ may wish to
consult \cite{St} or perhaps \cite{F3} or \cite{F5} which give worked examples in our notation.

	The contribution of this section is the not
obvious fact that cycle index methods can be applied to $P(n,q)$, a
maximal parabolic subgroup of $GL(n,q)$. It is also interesting that
the $n \rightarrow \infty$ limiting separable, cyclic, and semisimple
probabilities are all different from the corresponding limits in
$GL(n,q)$.

	The following notation will be used throughout this
section. $N(d,q)$ will denote the number of monic degree $d$
irreducible polynomials over $F_q$. $N^{'}(d,q)$ is defined by
$N^{'}(d,q)=N(d,q)$ for $d>1$ and $N'(1,q)=N(1,q)-2$. The following
elementary lemmas will be useful.

\begin{lemma}\label{allpoly} \[ \prod_{\phi}
(1-\frac{u^{deg(\phi)}}{q^{deg(\phi)}}) = 1-u \] \end{lemma}

\begin{proof} Expanding $\frac{1}{1-\frac{u^{deg(\phi)}}
{q^{deg(\phi)}}}$ as a geometric series and using unique factorization
in $F_q[x]$, one sees that the coefficient of $u^d$ in the reciprocal
of the left hand side is $\frac{1}{q^d}$ times the number of monic
polynomials of degree $d$, hence 1. Comparing with the reciprocal of
the right hand side completes the proof. \end{proof}

\begin{lemma} \label{Darboux} (Darboux \cite{O}) Suppose that $f(u)$
is analytic for $|u|<r,r>0$ and has only simple poles on
$|u|=r$. Letting $w_j$ denote the poles, and $g_j(u)$ be such that
$f(u)=\frac{g_j(u)}{1-u/w_j}$ and $g_j(u)$ is analytic near $w_j$, one
has that as $n \rightarrow \infty$, the coefficient of $u^n$ in $f(u)$ is \[ \sum_j
\frac{g_j(w_j)}{w_j^n}+o(1/r^n).\] \end{lemma}

\subsection{Separable Matrices}

	An element $\alpha$ in $A(n,q)$ is called separable if its
characteristic polynomial is square free. We remark that the results
of this subsection (but not those of the cyclic and semisimple
subsections) could also by obtained by using the cycle index for
$GL(n,q)$.

	Let $s_A(n,q)$ be the proportion of separable elements in
$A(n,q)$ and define the generating function $S_A(u,q) =
\sum_{n=0}^{\infty} u^n s_A(n,q)$ and let $S_{GL}(u,q)$ be the
corresponding generating function for $GL$.

\begin{theorem} \label{sepgen}
\begin{eqnarray*}
S_A(u,q) & = & \prod_{d \geq 1} (1+\frac{u^d}{q^d-1})^{N'(d,q)}\\
& = & \frac{S_{GL}(u,q)}{(1+\frac{u}{q-1})}
\end{eqnarray*}
\end{theorem}

\begin{proof} The first equality follows from Corollary \ref{cycleindex}
together with the fact that an element $\alpha$ of $GL(n,q)$ is separable if and only if all $\lambda_{\phi}(\alpha)$ have size at most
1. The second equality follows from the cycle index for
$GL(n,q)$. \end{proof}

	Wall \cite{W} shows that $S_{GL}(u,q)$ is analytic in $|u|<q$
except for a simple pole at $u=1$ which has residue $1-1/q$, implying
that $s_{GL}(\infty,q)=1-\frac{1}{q}$. To compute the corresponding
limit for the affine case note that for $q >2$, $S_A(u,q)$ is analytic in
$|u|<q-1$ except for a simple pole at $u=1$ which has residue
$\frac{1-\frac{1}{q}}{1+\frac{1}{q-1}}$. For $q=2$, since $N'(1,q)=0$
the function $S_A(u,q)$ is analytic in $|u|<q$ except for a simple
pole at $u=1$ which has residue
$\frac{1-\frac{1}{q}}{1+\frac{1}{q-1}}$. Thus
$s_A(\infty,q)=\frac{1-\frac{1}{q}}{1+\frac{1}{q-1}}$.

	Theorem \ref{boundsep} bounds the convergence rate of
$s_A(n,q)$ to its limit. The statement of the theorem is not fully
simplified so as to make the proof easier to follow.

\begin{theorem} \label{boundsep} Let $c=3/2$ and
$K^{+}=\frac{kc}{c-1}(1+\frac{q-2}{c^2}) (\frac{c}{q-c})$. Then 

\begin{enumerate}
\item $|s_A(n,q)-s_A(\infty,q)| \leq \frac{2K^{+}c (c/(q-1))^n}{(c-1)(1-c/(q-1))} + \frac{1}{(1-1/(q-1))(q-1)^{n+1}}$ for $q>2$.
\item $|s_A(n,q)-s_A(\infty,q)| \leq \frac{2K^{+}(c/q)^{n+1}}{1-(c/q)^2}$ for $q=2$.
\end{enumerate}
\end{theorem}

\begin{proof} Theorem \ref{sepgen} implies that

\[ (1-u) S_A(u,q) = \frac{(1-u)}{(1+\frac{u}{q-1})} S_{GL}(u,q).\]
Taking coefficients of $u^{m+1}$ on both sides and using the fact that $s_{GL}(1,q)=s_{GL}(0,q)$ one obtains that

\[ |s_A(m+1,q)-s_A(m,q)| \leq (\frac{1}{q-1})^m [\sum_{i=1}^m (q-1)^i
|s_{GL}(i+1,q)-s_{GL}(i,q)|]+\frac{1}{(q-1)^{m+1}}.\] For any $c,K^{+}$ as above, page 273
of \cite{W} gives the bound \[ |s_{GL}(i,q)-s_{GL}(\infty,q)| \leq
K^{+}(c/q)^i. \] Thus by the triangle inequality,
\[ |s_{GL}(i+1,q)-s_{GL}(i,q)| \leq 2K^+(c/q)^i.\] Now summing over $m
\geq n$ and using the triangle inequality gives

\begin{eqnarray*}
|s_A(n,q)-s_A(\infty,q)| & \leq &  2 \sum_{m \geq n} (\frac{1}{q-1})^m
[\sum_{i=1}^m (\frac{q-1}{q})^i K^{+} c^i] + \sum_{m \geq n} \frac{1}{(q-1)^{m+1}}\\
& \leq &  2 \sum_{m \geq n} (\frac{1}{q-1})^m
[\sum_{i=1}^m K^{+} c^i] + \frac{1}{(1-1/(q-1))(q-1)^{n+1}}\\
& \leq & \frac{2K^{+}c}{c-1} \sum_{m \geq n} (\frac{c}{q-1})^m + \frac{1}{(1-1/(q-1))(q-1)^{n+1}}\\
& = &  \frac{2K^{+}c (c/(q-1))^n}{(c-1)(1-c/(q-1))} + \frac{1}{(1-1/(q-1))(q-1)^{n+1}}.
\end{eqnarray*} For $q=2$ observe that taking coefficients of $u^{m+2}$ on both sides of
\[ (1-u^2) S_A(u,q) = (1-u) S_{GL}(u,q)\] gives that \[
|s_A(m+2,q)-s_A(m,q)| \leq |s_{GL}(m+2,q)-s_{GL}(m+1,q)|.\] Now use
Wall's bound for $|s_{GL}(i,q)-s_{GL}(\infty,q)|$ to conclude that

\begin{eqnarray*}
|s_A(n,q)-s_A(\infty,q)| & \leq & \sum_{m \geq n \atop m-n \ even} |s_A(m+2,q)-s_A(m,q)|\\
& \leq & \sum_{m \geq n \atop m-n \ even} |s_{GL}(m+2,q)-s_{GL}(m+1,q)|\\
& \leq & 2  \sum_{m \geq n \atop m-n \ even} K^+(c/q)^{m+1}\\
& \leq & \frac{2K^{+}(c/q)^{n+1}}{1-(c/q)^2}.
\end{eqnarray*}
\end{proof}

\subsection{Cyclic matrices}

	An element $\alpha$ of $GL(n,q)$ is called cyclic if its
characteristic polynomial is equal to its minimal polynomial. Let
$c_A(n,q)$ be the proportion of cyclic elements in $A(n,q)$ and let
$c_A(\infty,q)$ be the $n \rightarrow \infty$ limit of $c_A(n,q)$. Let
$C_A(u,q) = \sum_{n=0}^{\infty} u^n c_A(n,q)$ and let $C_{GL}(u,q)$ be
the corresponding generating function for $GL$.

\begin{theorem} \label{cycgen} 
\begin{eqnarray*}
C_A(u,q) & = & \frac{1}{1-u/q} \prod_{d \geq 1} (1+\frac{u^d}{(q^d-1)(1-u^d/q^d)})^{N'(d,q)}\\
& = & \frac{1}{1-u/q} \frac{1}{1+\frac{u}{(q-1)(1-u/q)}} C_{GL}(u,q).
\end{eqnarray*}
\end{theorem}

\begin{proof} The first equality follows from Corollary \ref{cycleindex}
together with the fact that an element $\alpha$ of $GL(n,q)$ is cyclic
if and only if all $\lambda_{\phi}(\alpha)$ have at most 1 part. The
second equality follows from the cycle index for $GL(n,q)$. \end{proof}

\begin{cor} \label{cyclim} \[ c_A(\infty,q) =
\frac{1-\frac{1}{q}}{1-\frac{1}{q}+\frac{1}{q^2}}
\frac{1-1/q^5}{1+1/q^3}.\] \end{cor}

\begin{proof} Wall \cite{W} shows that $C_{GL}(u,q)$ is analytic in
$|u|<q^2$ except for a simple pole at $u=1$ which has residue
$\frac{1-1/q^5}{1+1/q^3}$. Theorem \ref{cycgen} implies that
$C_{A}(u,q)$ is analytic in $|u|<q^2-q$ except for a simple pole at
$u=1$ which has residue
$\frac{1-\frac{1}{q}}{1-\frac{1}{q}+\frac{1}{q^2}}
\frac{1-1/q^5}{1+1/q^3}$. Thus $c_A(\infty,q)$ is equal to
$\frac{1-\frac{1}{q}}{1-\frac{1}{q}+\frac{1}{q^2}}
\frac{1-1/q^5}{1+1/q^3}$. \end{proof}

	For large $q$ the limit in Theorem \ref{cyclim} is of the
form $1-1/q^2+O(1/q^3)$. It would be interesting to understand this in
terms of algebraic geometry; the $1-1/q^3+O(q^4)$ behavior in the case
of $GL(n,q)$ is a finite analog of Steinberg's result that the variety
of regular semisimple elements has codimension three \cite{St}. See
\cite{NP3} for further discussion.

	Adapting a trick of Wall \cite{W} from the $GL$ case gives
instant bounds on the convergence rate of $c_A(n,q)$ to its limit. 

\begin{lemma} \label{trick} \[ (1-u) C_A(u,q) = (1-\frac{u}{q})
S_A(\frac{u}{q},q).\] \end{lemma}

\begin{proof} Theorems \ref{cycgen}, \ref{sepgen} and Lemma
\ref{allpoly} imply that
\begin{eqnarray*}
(1-u) C_A(u,q) & = & \frac{1}{1-u/q} \prod_{d \geq 1}
(1+\frac{u^d}{q^d(1-u^d/q^d)(1-1/q^d)})^{N'(q,d)} \prod_{d \geq 1}
(1-u^d/q^d)^{N(d,q)}\\
& = & (1-u/q) \prod_{d \geq 1} (1+\frac{u^d}{q^d(q^d-1)})^{N'(q,d)}\\
& = & (1-u/q) S_A(\frac{u}{q},q).
\end{eqnarray*}
\end{proof}

\begin{cor} \[ |c_A(\infty,q)-c_A(n,q)| \leq \frac{1}{q^{n+1}(1-1/q)}\]
\end{cor} 

\begin{proof} Taking coefficients of $u^{m+1}$ on both sides of the
equation in Lemma \ref{trick} gives that

\[ c_A(m+1,q)-c_A(m,q) = \frac{1}{q^{m+1}} [s_A(m+1,q)-s_A(m,q)].\]
Using the triangle inequality and the fact that $0 \leq
s_A(m+1,q),s_A(m,q) \leq 1$, it follows that

\begin{eqnarray*}
|c_A(\infty,q)-c_A(n,q)| & \leq & \sum_{m=n}^{\infty} |c_A(m+1,q)-c_A(m,q)|\\
& = & \sum_{m=n}^{\infty} \frac{1}{q^{m+1}} |s_A(m+1,q)-s_A(m,q)|\\
& \leq & \sum_{m=n}^{\infty} \frac{1}{q^{m+1}}\\
& = & \frac{1}{q^{n+1}(1-1/q)}.
\end{eqnarray*}
\end{proof} 

\subsection{Semisimple matrices}

	An element $\alpha$ of $GL(n,q)$ is called semisimple if it
diagonalizable over the algebraic closure of $F_q$. Treatments of
semisimple probabilities without generating functions appear in
\cite{IsKanSp} and \cite{GL}, which proves the lovely result that if
$G$ is a simple Chevalley group, then the probability of not being
semisimple is at most $3/(q-1)+2/(q-1)^2$.

	 A crude asymptotic understanding of the behavior of the
proportion of semisimple elements in $GL(n,q)$ is in \cite{St}, who
used generating functions. The paper \cite{F5} uses one of the
Rogers-Ramanujan identities to show that $n \rightarrow \infty$
probability that an element of $GL(n,q)$ is semisimple is \[
\prod_{r=1 \atop r=0,\pm 2 (mod \ 5)}^{\infty}
\frac{(1-\frac{1}{q^{r-1}})}{(1-\frac{1}{q^r})}. \] The paper
\cite{FNP} gives bounds for finite $n$.

	Let $ss_A(n,q)$ be the probability that an element of $A(n,q)$
is semisimple and let $ss_A(\infty,q)$ be the $n \rightarrow \infty$
limit of $ss_A(n,q)$. Let $SS_A(u,q) = \sum_{n=0}^{\infty} u^n
ss_A(n,q)$ and let $SS_{GL}$ be the corresponding generating function
for $GL$.

\begin{theorem} \label{genss} 
\begin{eqnarray*}
SS_A(u,q) & = & (\sum_{k \geq 0} \frac{u^k}{q^{k^2+k}
(\frac{1}{q})_k}) \prod_{d \geq 1} (\sum_{k \geq 0} \frac{u^{kd}}{q^{kd}
(\frac{1}{q^d})_k})^{N'(d,q)}\\
& = & \frac{(\sum_{k \geq 0} \frac{u^k}{q^{k^2+k}
(\frac{1}{q})_k})}{(\sum_{k \geq 0} \frac{u^k}{q^{k^2}
(\frac{1}{q})_k})} SS_{GL}(u,q).
\end{eqnarray*}
\end{theorem} 

\begin{proof} An element $\alpha$ of $GL(n,q)$ is semisimple if and only if all $\lambda_{\phi}(\alpha)$ have at most one column. Now use Corollary \ref{cycleindex}. \end{proof}

	To calculate $ss_A(\infty,q)$ we will use (with $q$ replaced
by $1/q$) the well known Rogers-Ramanujan identities (see \cite{A} for
discussion)

\[ 1+\sum_{n=1}^{\infty} \frac{q^{n^2}}{(1-q)(1-q^2) \cdots (1-q^n)} =
\prod_{n=1}^{\infty} \frac{1}{(1-q^{5n-1})(1-q^{5n-4})} \]

\[ 1+\sum_{n=1}^{\infty} \frac{q^{n(n+1)}}{(1-q)(1-q^2) \cdots
(1-q^n)} = \prod_{n=1}^{\infty} \frac{1}{(1-q^{5n-2})(1-q^{5n-3})}. \]
  
\begin{cor} \[ ss_A(\infty,q) = \frac{\prod_{r=1 \atop r=0,\pm 1 \ mod \ 5}^{\infty} (1-1/q^r) \prod_{r=1 \atop r=0,\pm 2 \ mod \ 5}^{\infty} (1-1/q^{r-1})}{\prod_{r=1 \atop r=0,\pm 2 \ mod \ 5}^{\infty} (1-1/q^{r})^2} \] \end{cor}

\begin{proof} The discussion in \cite{FNP} shows that $SS_{GL}(u,q)$
is analytic within a circle of radius greater than 1, except for a
simple pole at $u=1$. Since \[ \frac{(\sum_{k \geq 0}
\frac{u^k}{q^{k^2+k} (\frac{1}{q})_k})}{(\sum_{k \geq 0}
\frac{u^k}{q^{k^2} (\frac{1}{q})_k})} \] is also analytic within a
circle of radius greater than 1, it follows that \[ ss_A(\infty,q) =
\frac{(\sum_{k \geq 0} \frac{1}{q^{k^2+k} (\frac{1}{q})_k})}{(\sum_{k
\geq 0} \frac{1}{q^{k^2} (\frac{1}{q})_k})} ss_{GL}(\infty,q).\] Now
simply use both Rogers-Ramanujan identities and the formula for
$ss_{GL}(\infty,q)$ in \cite{F5} stated at the beginning of this
subsection. \end{proof}

	Bounding the convergence rate of $ss_{A}(n,q)$ through
generating functions is an involved analytic excursion which we
omit. The following elementary bound is sufficient for practical
purposes, given the results of the previous subsections and the fact
that $ss_A(n,q)$ can be computed explicitly for small $n$ from the
generating function.

\begin{theorem} \[ s_A(n,q) \leq ss_A(n,q) \leq s_A(n,q) +
(1-c_A(n,q)) .\] \end{theorem}

\begin{proof} The first inequality follows because semisimple matrices
are separable. The second inequality follows because a matrix which is
not separable is either not cyclic or not semisimple. \end{proof}

\section*{Acknowledgements} The author thanks Persi Diaconis, Robert Guralnick, Peter M. Neumann, and Cheryl E. Praeger for conversations related to this
paper. This research was done at Stanford University with the support of an NSF
Postdoctoral Fellowship.

\end{document}